\documentstyle[11pt,twoside]{article}
\include{graphicx}
\title{A note on a class of exact solutions for a doubly anharmonic oscillator}
\author{\sc R. B.\ Paris \\
{\em Division of Computing and Mathematics}, \\
{\em University of Abertay Dundee, Dundee DD1 1HG, UK}
}
\begin{document}
\def\f#1#2{\mbox{${\textstyle \frac{#1}{#2}}$}}
\def\dfrac#1#2{\displaystyle{\frac{#1}{#2}}}
\def\boldal{\mbox{\boldmath $\alpha$}}
{\newcommand{\Sgoth}{S\;\!\!\!\!\!/}
\newcommand{\bee}{\begin{equation}}
\newcommand{\ee}{\end{equation}}
\newcommand{\lam}{\lambda}
\newcommand{\ka}{\kappa}
\newcommand{\al}{\alpha}
\newcommand{\fr}{\frac{1}{2}}
\newcommand{\fs}{\f{1}{2}}
\newcommand{\g}{\Gamma}
\newcommand{\br}{\biggr}
\newcommand{\bl}{\biggl}
\newcommand{\ra}{\rightarrow}
\newcommand{\mbint}{\frac{1}{2\pi i}\int_{c-\infty i}^{c+\infty i}}
\newcommand{\mbcint}{\frac{1}{2\pi i}\int_C}
\newcommand{\mboint}{\frac{1}{2\pi i}\int_{-\infty i}^{\infty i}}
\newcommand{\gtwid}{\raisebox{-.8ex}{\mbox{$\stackrel{\textstyle >}{\sim}$}}}
\newcommand{\ltwid}{\raisebox{-.8ex}{\mbox{$\stackrel{\textstyle <}{\sim}$}}}
\renewcommand{\topfraction}{0.9}
\renewcommand{\bottomfraction}{0.9}
\renewcommand{\textfraction}{0.05}
\newcommand{\mcol}{\multicolumn}
\date{}
\maketitle
\pagestyle{myheadings}
\markboth{\hfill \sc R. B.\ Paris  \hfill}
{\hfill \sc Doubly anharmonic oscillator\hfill}
\begin{abstract}
We examine a class of exact solutions for the eigenvalues and eigenfunctions of a doubly anharmonic oscillator defined by the potential $V(x)=\omega^2/2 x^2+\lambda x^4/4+\eta x^6/6$, $\eta>0$. These solutions hold provided  certain constraints on the coupling parameters $\omega^2$, $\lambda$ and $\eta$ are satisfied.
\vspace{0.4cm}

\noindent {\bf Mathematics Subject Classification:} 34A05, 81Q05
\vspace{0.3cm}

\noindent {\bf Keywords:}  Schr\"odinger equation, doubly anharmonic potential, exact solutions
\end{abstract}

\vspace{0.3cm}

\noindent $\,$\hrulefill $\,$

\vspace{0.2cm}

\begin{center}
{\bf 1. \  Introduction}
\end{center}
\setcounter{section}{1}
\setcounter{equation}{0}
\renewcommand{\theequation}{\arabic{section}.\arabic{equation}}
The anharmonic oscillator we consider in this note has the potential
$V(x)=\fs\omega^2x^2+\f{1}{4}\lambda x^4+\f{1}{6}\eta x^6$,
where $\omega^2$, $\lambda$ and $\eta>0$ are real parameters.
With $E$ denoting the energy eigenvalue, the associated Schr\"odinger equation is
\bee\label{e11}
\frac{d^2\psi(x)}{dx^2}+(2E-\omega^2x^2-\fs\lambda x^4-\f{1}{3}\eta x^6)\,\psi(x)=0\qquad (-\infty<x<\infty)
\ee
subject to the boundary conditions
\[\lim_{x\ra\pm\infty}\psi(x)=0.\]
In \cite{F1}, Flessas obtained exact solutions for the energy eigenvalues and eigenfunctions (for the ground state and first excited state) when certain constraints on the coupling parameters are satisfied. Singh {\it et al.\/} \cite{S}
subsequently showed how an infinite set of such solutions could be constructed, although they presented details only of the eigenfunctions with 2 and 3 nodes. Again, these solutions are obtained when a sequence of constraints on the coupling parameters is satisfied; a finite number of such solutions is engendered by
the fulfillment of each set of constraints.

In this note we consider the case of eigensolutions with up to $2N$ and $2N+1$ nodes, where $N\leq 3$, in more detail.
\vspace{0.6cm}

\begin{center}
{\bf 2. \  Solution scheme}
\end{center}
\setcounter{section}{2}
\setcounter{equation}{0}
\renewcommand{\theequation}{\arabic{section}.\arabic{equation}}
We define the quantities
\bee\label{e20}
a=\frac{\lambda}{4}\bl(\frac{3}{\eta}\br)^\fr,\qquad b=\bl(\frac{\eta}{3}\br)^\fr,\qquad c=\omega^2+(3\eta)^\fr-\frac{3\lambda^2}{16\eta}.
\ee
With the substitution
\bee\label{e12}
\psi(x)=\exp\,[-\fs ax^2-\f{1}{4}bx^4]\,y(x)
\ee
in (\ref{e11}), we obtain the differential equation
\bee\label{e13}
y''(x)-2(ax+bx^3) y'(x)+(2E-a-cx^2) y(x)=0.
\ee

We look for polynomial solutions of (\ref{e13}) of the form
\bee\label{e14}
y(x)=\sum_{n=0}^N A_n x^{2n+\epsilon}\qquad (A_0=1;\ N=0, 1, 2, \ldots),
\ee
where $\epsilon=0$ (even solutions) or $\epsilon=1$ (odd solutions) and the $A_n$ are coefficients to be determined. Substitution of (\ref{e14}) into the differential equation (\ref{e13}) then yields
\[\sum_{n=0}^N A_n\bl\{(2n+\epsilon)(2n-1+\epsilon) x^{2n-2}+(2E-a-2a(2n+\epsilon) x^{2n}\]
\[-(c+2b(2n+\epsilon))x^{2n+2}\br\}=0.\]
The requirement that the constant terms and the coefficient of $x^{2N+2}$ should vanish produces
\[(1+\epsilon)(2+\epsilon)A_1+2E-a(1+2\epsilon)=0,\]
\bee\label{e15a}
c+2b(2N+\epsilon)=0.
\ee

Thus we have the conditions
\bee\label{e15}
E=\fs a(1+2\epsilon)-\fs(1+\epsilon)(2+\epsilon) A_1
\ee
and
\bee\label{e16}
\gamma:=\bl(\frac{3}{\eta}\br)^\fr \bl(\frac{3\lambda^2}{16\eta}-\omega^2\br)=4N+3+2\epsilon,
\ee
together with the $N$ equations for the coefficients $A_n$
\[(2n+1+\epsilon)(2n+2+\epsilon) A_{n+1}+(2E-a-2a(2n+\epsilon)) A_n\hspace{3cm}\]
\bee\label{e17}
\hspace{4cm}-(c+2b(2n-2+\epsilon)) A_{n-1}=0\qquad (1\leq n\leq N),
\ee
with $A_n=0$ for $n>N$. 

In what follows, we label the eigenvalues associated with the even solutions by $E_{2m}$, $m=0, 1, 2, \ldots$, where $E_0$ denotes the ground state eigenvalue, and those associated with the odd solutions by $E_{2m+1}$. A similar indexing applies to the eigenfunctions $\psi(x)$, which will possess either $2m$ or $2m+1$ nodes on the interval $(-\infty, \infty)$, respectively.
\vspace{0.5cm}

\noindent (i) {\em The case $N=0$}.\ \ \ When $N=0$, we obtain from (\ref{e15}) and (\ref{e16}) the values
\[E_0=\fs a,\ \ \gamma=3\qquad \mbox{and}\qquad E_1=\f{3}{2}a,\ \ \gamma=5.\]
The associated eigenfunctions are
\[\psi_0(x)=\exp\,[-\fs ax^2-\f{1}{4}bx^4]\qquad\mbox{and}\qquad \psi_1(x)=x\,\exp\,[-\fs ax^2-\f{1}{4}bx^4],\]
which possess 0 and 1 node, respectively.
These are the special solutions of (\ref{e11}) obtained by Flessas \cite{F1}.
\vspace{0.5cm}

\noindent (ii) {\em The case $N=1$}.\ \ \ When $N=1$, we obtain from (\ref{e15a}), (\ref{e15}) and (\ref{e16})  for $\epsilon=0$ the values
\[E=\fs a-A_1,\quad \gamma=7\ \ \ (c=-4b)\]
and from (\ref{e17}) the single equation
\[(2E-5a)A_1-c=0.\]
Substitution of the above values of $E$ and $c$ yields the quadratic
\[A_1^2+2a A_1+2b=0,\]
with solutions
\[A_1=-a\pm\sqrt{a^2+2b}.\]

The positive root is associated with an eigenfunction with no node and so is a ground-state value \cite[p.~21]{F2}, viz.
\[E_0=\f{3}{2}a-\sqrt{a^2+2b},\quad \psi_0(x)=\{1+(\sqrt{a^2+2b}-a)x^2\} \exp\,[-\fs ax^2-\f{1}{4}bx^4].\]
The negative root corresponds to an eigenfunction with 2 nodes and represents a case of the first even excited state:
\[E_2=\f{3}{2}a+\sqrt{a^2+2b},\quad \psi_2(x)=\{1-(\sqrt{a^2+2b}+a)x^2\} \exp\,[-\fs ax^2-\f{1}{4}bx^4].\]
Both cases are associated with the value $\gamma=7$, which gives a constraint between the three parameters $\omega^2$, $\lambda$ and $\eta$.

When $\epsilon=1$, a similar procedure shows that 
\[E=\f{3}{2}a-3A_1,\quad \gamma=9\]
and
\[A_1=-\f{1}{3}a\pm\f{1}{3} \sqrt{a^2+6b}.\]
Again, the positive root corresponds to the lowest odd eigenfunction (with 1 node) and the negative root to the second odd eigenfunction (with 3 nodes). Thus, we have
\[E_1=\f{5}{2}a-\sqrt{a^2+6b},\quad\psi_1(x)=x\{1+\f{1}{3}(\sqrt{a^2+6b}-a) x^2\} \exp\,[-\fs ax^2-\f{1}{4}bx^4],\]
\[E_3=\f{5}{2}a+\sqrt{a^2+6b},\quad\psi_3(x)=x\{1-\f{1}{3}(\sqrt{a^2+6b}+a) x^2\} \exp\,[-\fs ax^2-\f{1}{4}bx^4].\]
These special solutions of (\ref{e11}) were given in \cite{S}.
\vspace{0.6cm}

\begin{center}
{\bf 3. \  Exact solutions for the case $N=2$}
\end{center}
\setcounter{section}{3}
\setcounter{equation}{0}
\renewcommand{\theequation}{\arabic{section}.\arabic{equation}}
\noindent When $N=2$, we have from (\ref{e15a}), (\ref{e15}) and (\ref{e16}) for even modes ($\epsilon=0$)
\[E=\fs a-A_1,\quad \gamma=11\ \ \ (c=-8b)\]
together with, from (\ref{e17}), the two equations for the coefficients
\[12A_2+(2E-5a) A_1-c=0,\qquad (2E-9a) A_2-(c+4b) A_1=0.\]
Substitution of the above values of $E$ and $c$ then leads after some straightforward algebra to
\bee\label{e18}
A_1^3+6a A_1^2+8(a^2-2b) A_1-16ab=0,\qquad A_2=\frac{2bA_1}{A_1+4a}.
\ee

The cubic equation for $A_1$ can be written in its reduced form as
\[\chi^3+p\chi+q=0,\qquad A_1=\chi-2a,\]
where $p=-4(a^2+4b)$ and $q=16ab$. The discriminant
\[\Delta=-4p^3-27q^2=256\{(a^2+4b)^3-27a^2b^2\}>0,\]
so that (\ref{e18}) has three real roots. Moreover, inspection of the coefficients of (\ref{e18}) shows that there are two negative roots and one positive root. With
\[P=(-\f{4}{3}p)^\fr=4\bl(\frac{a^2+4b}{3}\br)^{\!\fr},\quad Q=\bl(\frac{-27q^2}{4p^3}\br)^{\!\fr}=ab \bl(\frac{a^2+4b}{3}\br)^{\!\!-\frac{3}{2}},\]
\bee\label{edef}
\theta=\f{1}{3} \arcsin Q,
\ee
the roots are then given by
\[\chi_k=P \sin\,(\theta+\f{2}{3}\pi k) \qquad (k=0, 1, 2).\]
It is easy to establish that the positive value of $A_1$ corresponds to the root $\chi_1$.

Then we have the following even exact solutions: ({\it a})\ a ground-state eigenvalue and eigenfunction given by
\bee\label{e21}
E_0=\f{5}{2}a-P \sin\,(\theta+\f{2}{3}\pi),\quad \psi_0(x)=\{1+A_1x^2+A_2x^4\} \exp\,[-\fs ax^2-\f{1}{4}bx^4],
\ee
where
\[A_1=P \sin\,(\theta+\f{2}{3}\pi)-2a>0,\qquad A_2=\frac{2b(P \sin\,(\theta+\f{2}{3}\pi)-2a)}{P \sin\,(\theta+\f{2}{3}\pi)+2a}>0;\]
({\it b})\ a first even excited state (2 nodes) with eigenvalue and eigenfunction given by
\bee\label{e22}
E_2=\f{5}{2}a-P \sin\,\theta,\quad \psi_2(x)=\{1-|A_1|x^2-|A_2|x^4\} \exp\,[-\fs ax^2-\f{1}{4}bx^4],
\ee
where
\[A_1=P \sin\,\theta-2a<0,\qquad A_2=\frac{2b(P \sin\,\theta-2a)}{P \sin\,\theta+2a}<0;\]
and ({\it c})\ a second even excited state (4 nodes) with eigenvalue and eigenfunction given by
\bee\label{e23}
E_4=\f{5}{2}a-P \sin\,(\theta+\f{4}{3}\pi),\quad \psi_4(x)=\{1-|A_1|x^2+A_2x^4\} \exp\,[-\fs ax^2-\f{1}{4}bx^4],
\ee
where\footnote{The quantity $P \sin\,(\theta+\f{4}{3}\pi)+2a<0$, since $P |\sin (\theta+\frac{4}{3}\pi)|>(4a/\surd 3) |\sin \frac{4}{3}\pi|=2a$ when $0\leq\theta\leq\frac{\pi}{6}$.}
\[A_1=P \sin\,(\theta+\f{4}{3}\pi)-2a<0,\qquad A_2=\frac{2b(P \sin\,(\theta+\f{4}{3}\pi)-2a)}{P \sin\,(\theta+\f{4}{3}\pi)+2a}>0.\]

A similar treatment for odd eigensolutions ($\epsilon=1$) yields
\[E=\f{3}{2} a-3A_1,\quad \gamma=13\ \ \ (c=-10b)\]
and
\[(3A_1)^3+6a (3A_1)^2+8(a^2-4b) (3A_1)-48ab=0,\qquad A_2=\frac{2bA_1}{3A_1+4a}.\]
The reduced cubic is
\[\chi^3+p\chi+q=0,\qquad 3A_1=\chi-2a,\]
where $p=-4(a^2+8b)$ and $q=16ab$. The three real roots are given by
\[\chi_k=P \sin\,(\theta+\f{2}{3}\pi k)\qquad (k=0, 1, 2),\]
where now
\[P=4\bl(\frac{a^2+8b}{3}\br)^{\!\fr},\quad Q=ab \bl(\frac{a^2+8b}{3}\br)^{\!\!-\frac{3}{2}},\quad\theta=\f{1}{3} \arcsin Q,\]

Then we have the following odd exact solutions: ({\it a})\ a first excited eigenvalue and eigenfunction (with 1 node) given by
\bee\label{e24}
E_1=\f{7}{2}a-P \sin\,(\theta+\f{2}{3}\pi),\quad \psi_1(x)=x\{1+A_1x^2+A_2x^4\} \exp\,[-\fs ax^2-\f{1}{4}bx^4],
\ee
where
\[A_1=\frac{P}{3} \sin\,(\theta+\f{2}{3}\pi)-\frac{2a}{3}>0,\qquad A_2=\frac{2b(P \sin\,(\theta+\f{2}{3}\pi)-2a)}{3(P \sin\,(\theta+\f{2}{3}\pi)+2a)}>0;\]
({\it b})\ a second excited state (3 nodes) with eigenvalue and eigenfunction given by
\bee\label{e25}
E_3=\f{7}{2}a-P \sin\,\theta,\quad \psi_3(x)=x\{1-|A_1|x^2-|A_2|x^4\} \exp\,[-\fs ax^2-\f{1}{4}bx^4],
\ee
where
\[A_1=\frac{P}{3} \sin\,\theta-\frac{2a}{3}<0,\qquad A_2=\frac{2b(P \sin\,\theta-2a)}{3(P \sin\,\theta+2a)}<0;\]
and ({\it c})\ a third excited state (5 nodes) with eigenvalue and eigenfunction given by
\bee\label{e26}
E_5=\f{7}{2}a-P \sin\,(\theta+\f{4}{3}\pi),\quad \psi_5(x)=x\{1-|A_1|x^2+A_2x^4\} \exp\,[-\fs ax^2-\f{1}{4}bx^4],
\ee
where
\[A_1=\frac{P}{3} \sin\,(\theta+\f{4}{3}\pi)-\frac{2a}{3}<0,\qquad A_2=\frac{2b(P \sin\,(\theta+\f{4}{3}\pi)-2a)}{3(P \sin\,(\theta+\f{4}{3}\pi)+2a)}>0.\]
\vspace{0.6cm}

\begin{center}
{\bf 4. \  Exact solutions for the case $N=3$}
\end{center}
\setcounter{section}{4}
\setcounter{equation}{0}
\renewcommand{\theequation}{\arabic{section}.\arabic{equation}}
When $N=3$, we have from (\ref{e15a}), (\ref{e15}) and (\ref{e16}) for even modes ($\epsilon=0$)
\[E=\fs a-A_1,\quad \gamma=15\ \ \ (c=-12b)\]
and from (\ref{e17}) the equations for the coefficients
\begin{eqnarray*}
12A_2+(2E-5a)A_1-c&=&0\\
30A_3+(2E-9a)A_2-(c+4b)A_1&=&0\\
(2E-13a)A_3-(c+8b)A_2&=&0.
\end{eqnarray*}
This produces the quartic equation for $w=A_1$ given by
\bee\label{e27}
w^4+12aw^3+4(11a^2-15b)w^2+24a(2a^2-11b)w-36b(4a^2-5b)=0,
\ee
\[A_2=\frac{4bA_1(A_1+6a)}{(A_1+4a)(A_1+6a)-30b},\qquad A_3=\frac{2bA_2}{A_1+6a}~.\]

The reduced quartic is
\bee\label{e2q}
\chi^4+p\chi^2+q\chi+r=0,\qquad A_1=\chi-3a,
\ee
where $p=-10(a^2+6b)$, $q=96ab$ and $r=9(a^4+12a^2b+20b^2)$. This equation possesses four real\footnote{This follows from the fact that $p<0$, $D=64r-16p^2=-1024(a^4+12a^2b+45b^2)<0$ and the discriminant $\Delta$ defined in \cite[Eq.~(1.11.17)]{DLMF} is 
$\Delta=3^2\times 2^{16} (a^{12}+36a^{10}b+402a^8b^2+1848a^6b^3+3897a^4b^4+35100a^2b^5+40500b^6)>0$.} roots.
Although it is possible to express the roots in algebraic form, it was found that the resulting expressions were too complicated to be of practical use. In this case, we shall content ourselves with a numerical solution of the quartic equation (\ref{e27}).

If we choose, for example, $\lambda=0.50$, $\eta=0.03$ ($a=1.25$, $b=0.10$) then, from the constraint $\gamma=13$, we have $\omega^2=0.0625$. The largest root of (\ref{e27}) together with the corresponding values of $A_2$ and $A_3$ are
\[A_1=0.264080,\quad A_2=0.021656,\quad A_3=0.000558.\]
We note that all the coefficients are positive and so this will result in a ground-state eigenfunction (no node) with the eigenvalue
$E_0=\f{5}{8}-A_1=0.360920$. 
Values of the other solutions of (\ref{e27}), which we label $A_n^{(2m)}$ ($0\leq m\leq 3$), and the associated even eigenvalues $E_{2m}$ are presented in Table 1
and the corresponding eigenfunctions are
\[\psi_{2m}(x)=\{1+\sum_{n=1}^3A_n^{(2m)}x^{2n}\} \exp\,[-\f{5}{8}x^2-\f{1}{40}x^4]\qquad(0\leq m\leq 3).\]
\begin{table}[t]
\caption{\footnotesize{Values of the coefficients when $\lambda=0.50$, $\eta=0.03$ and $\omega^2=0.0625$ obtained from (\ref{e27}) and the corresponding even eigenvalues $E_{2m}$. }}
\begin{center}
\begin{tabular}{l|l|l|l||l}
\hline
&&&&\\[-0.3cm]
\mcol{1}{c|}{$m$} & \mcol{1}{c|}{$A_1^{(2m)}$} & \mcol{1}{c|}{$A_2^{(2m)}$} & \mcol{1}{c||}{$A_3^{(2m)}$}  & \mcol{1}{c}{$E_{2m}$} \\ 
[.1cm]\hline
&&&&\\[-0.3cm]
0 & $+0.264080$ & $+0.021656$ & $+0.000558$ & $0.360920$\\
1 & $-1.887128$ & $-0.292761$ & $-0.010432$ & $2.512128$\\
2 & $-4.899957$ & $+1.859948$ & $+0.143071$ & $5.524957$\\
3 & $-8.476994$ & $+8.344491$ & $-1.708197$ & $9.101994$\\
[.2cm]\hline
\end{tabular}
\end{center}
\end{table}

For odd modes with $\epsilon=1$ we have
\[E=\f{3}{2}a-3A_1,\quad \gamma=17\ \ \ (c=-14b)\]
and the equations for the coefficients
\begin{eqnarray*}
20A_2+(2E-7a)A_1-(c+2b)&=&0\\
42A_3+(2E-11a)A_2-(c+6b)A_1&=&0\\
(2E-15a)A_3-(c+10b)A_2&=&0.
\end{eqnarray*}
This produces the quartic equation for $w=3A_1$ given by
\bee\label{e28}
w^4+12aw^3+(44a^2-100b)w^2+24a(2a^2-21b)w-108b(4a^2-7b)=0,
\ee
\[A_2=\frac{4bA_1(3A_1+6a)}{(3A_1+4a)(3A_1+6a)-42b},\qquad A_3=\frac{2bA_2}{3A_1+6a}~.\]
This equation also has four real roots and, for $\lambda=0.50$, $\eta=0.03$ and $\omega^2=0.0625$, the four sets of values of the coefficients $A_n^{(2m+1)}$ ($0\leq m\leq 3$) and the corresponding odd eigenvalues $E_{2m+1}$ are presented in Table 2. The associated eigenfunctions are
\[\psi_{2m+1}(x)=x\{1+\sum_{n=1}^3A_n^{(2m+1)}x^{2n}\} \exp\,[-\f{5}{8}x^2-\f{1}{40}x^4] \qquad (0\leq m\leq 3).\]
\begin{table}[t]
\caption{\footnotesize{Values of the coefficients when $\lambda=0.50$, $\eta=0.03$ and $\omega^2=0.0625$ obtained from (\ref{e28}) and the corresponding odd eigenvalues $E_{2m+1}$.}}
\begin{center}
\begin{tabular}{l|l|l|l||l}
\hline
&&&&\\[-0.3cm]
\mcol{1}{c|}{$m$} & \mcol{1}{c|}{$A_1^{(2m+1)}$} & \mcol{1}{c|}{$A_2^{(2m+1)}$} & \mcol{1}{c||}{$A_3^{(2m+1)}$}  & \mcol{1}{c}{$E_{2m+1}$} \\ 
[.1cm]\hline
&&&&\\[-0.3cm]
0 & $+0.243487$ & $+0.018657$ & $+0.000453$ & $1.144540$\\
1 & $-0.611015$ & $-0.100752$ & $-0.003556$ & $3.708044$\\
2 & $-1.690968$ & $+0.375069$ & $+0.030907$ & $6.947903$\\
3 & $-2.941504$ & $+1.800358$ & $-0.271852$ & $10.699513$\\
[.2cm]\hline
\end{tabular}
\end{center}
\end{table}

\vspace{0.4cm}

\begin{center}
{\bf 5. \  Summary}
\end{center}
\setcounter{section}{5}
\setcounter{equation}{0}
\renewcommand{\theequation}{\arabic{section}.\arabic{equation}}
We have examined in detail the cases $N=2$ and $N=3$ of exact solutions of the form
\bee\label{e31}
\psi(x)=\sum_{n=1}^N A_n x^{2n+\epsilon} \exp\,[-\fs a x^2-\f{1}{4}bx^4]\qquad (\epsilon=0, 1)
\ee
of the Schr\"odinger equation (\ref{e11}), where $a=\lambda (3/\eta)^\fr/4$ and $b=(\eta/3)^\fr$ with $\eta>0$. The cases $N=0$ and $N=1$ have been given earlier in \cite{F1} and \cite{S}, respectively. It is found that the solution (\ref{e31}) can only exist if a constraint on the coupling parameters $\omega^2$, $\lambda$ and $\eta$ is satisfied, viz.
\[\gamma=4N+3+2\epsilon=\bl(\frac{3}{\eta}\br)^{\!\fr} \bl(\frac{3\lambda^2}{16\eta}-\omega^2\br).\]
For a given $N$ and parity $\epsilon$ of the eigenfunction, two parameters are free to be chosen with the third then fixed by the above constraint. 

For each value of $N$ considered it is  found that $N+1$ eigenstates are produced with eigenvalues $E_0, E_2, \ldots , E_{2N}$ in the case of even modes and $E_1, E_3, \ldots , E_{2N+1}$ in the case of odd modes.
We present below a summary of the ground-state eigenvalues and eigenfunctions of type (\ref{e31}) with $\epsilon=0$ and normalised such that $A_0=1$:
\[N=0:\ \ \ E_0=\fs a,\ \ \gamma=3,\]
\[N=1:\ \ \ E_0=\f{3}{2}a-\sqrt{a^2+2b},\ \ \gamma=7; \qquad A_1=\sqrt{a^2+2b}-a, \]
\[N=2:\ \ \ E_0=\f{5}{2}a-P \sin (\theta+\f{2}{3}\pi),\ \ \ \gamma=11; \qquad A_1=P \sin (\theta+\f{2}{3}\pi)-2a,\ \ A_2=\frac{2bA_1}{A_1+4a},\]
\[N=3:\ \ \ E_0=\f{7}{2}a-\chi^*,\ \ \gamma=15,\]
where $P$ and $\theta$ are defined in (\ref{edef}) and $\chi^*$ denotes the largest root of the quartic (\ref{e2q}).
In the case $N=4$ only numerical solutions for the coefficients $A_n$ $(1\leq n\leq 3)$ are obtained for a specific choice of parameters.


\vspace{0.6cm}


\begin{thebibliography}{99}
\footnotesize{
\bibitem{F1}
G. P. Flessas, Exact solutions for a doubly anharmonic oscillator, Phys. Lett. {\bf 72A} (1979) 289--290.

\bibitem{F2}
G. P. Flessas, On the three-dimensional anharmonic oscillator, Phys. Lett. {\bf 78A} (1980) 19--21.

\bibitem{DLMF}
F. W. J. Olver, D. W. Lozier, R. F. Boisvert and C. W. Clark (eds.),    
{\it NIST Handbook of Mathematical Functions}, Cambridge University Press, Cambridge, 2010.

\bibitem{S}
V. Singh, A. Rampal, S. N. Biswas and K. Datta, A class of exact solutions for doubly anharmonic oscillators, Lett. Math. Phys. {\bf 4} (1980) 131--134.
}
\end{thebibliography}
\end{document}